\magnification=\magstep1 \hfuzz=1pt \tolerance=1000
\input amssym.def
\input amssym.tex
\font \normal=cmr10 scaled \magstep0 \font \mittel=cmr10 scaled
\magstep1 \font \gross=cmr10 scaled \magstep2
\def\Naturalsign{{\rm l\kern-.23em N}}
\def\Supp{\mathop {\rm Supp}\nolimits }
\def\Ass{\mathop {\rm Ass}\nolimits }
\def\Hom{\mathop {\rm Hom}\nolimits }
\def\Ext{\mathop {\rm Ext}\nolimits }
\def\injdim{\mathop {\rm injdim}}
\def\Spec{\mathop {\rm Spec}}
\def\depth{\mathop {\rm depth}}
\def\height{\mathop {\rm height}}
\baselineskip=15pt \gross \noindent On the set of associated
primes of a local cohomology module \normal
\bigskip
\bigskip
\mittel
\noindent
M. Hellus
\normal
\bigskip
\noindent Mathematische Fakult\"at der Universit\"at Regensburg

\noindent Address:    Universit\"at Regensburg, 93040 Regensburg, Germany

\noindent E-mail-address: michael.hellus@mathematik.uni-regensburg.de
\bigskip
\hrule
\medskip
{{\tt ABSTRACT:} Assume $R$ is a local Cohen-Macaulay
ring. It is shown that $\Ass_R (H^l_I(R))$ is finite for any ideal $I$ and
any integer $l$ provided $\Ass_R (H^2_{(x,y)}(R))$ is finite for any
$x,y\in R$ and $\Ass_R (H^3_{(x_1,x_2,y)}(R))$ is finite for any $y\in R$
and any regular sequence $x_1,x_2\in R$. Furthermore it is shown that
$\Ass_R (H^l_I(R))$ is always finite if $\dim (R)\leq 3$. The same
statement is even true for $\dim (R)\leq 4$ if $R$ is almost factorial.}
\medskip
\hrule
\bigskip
Cohomology theory is an important part of algebraic geometry.
If one considers local cohomology on an affine scheme with support
in a closed subset, everything can be expressed in terms of rings,
ideals and modules. More precisely, let $R$ be a noetherian ring and
$I$ an ideal of $R$ (determining a closed subset of $\Spec (R)$): In
this situation one studies the local cohomology modules $H^l_I(M)$,
where $l$ is a natural number and $M$ is any $R$-module. As these
local cohomology modules behave well under localisation, one
often restricts the above situation to the case $R$ is a local ring.

As the structure of local cohomology modules in general seems to
be quite mysterious, one tries to establish finiteness properties
providing a better understanding of these modules. Finiteness
properties of local cohomology modules have been studied by
several authors, see for example Brodmann/Lashgari Faghani [1],
Huneke/Koh [5],
Huneke/Sharp [6], Lyubeznik [8] and Singh [11]. For a survey
of results see Huneke [7].

Throughout this paper $(R,\goth m)$ is a local noetherian ring and
$I$ an ideal of $R$. We deal with the question, whether the set of
associated primes of every local cohomology module $H^l_I(R)$ is finite.
As local cohomology modules in general are not finitely generated, this
is an interesting question. For example if $R$ is a regular local
ring containing a field then $H^l_I(R)$ (for $l\geq 1$) is finitely
generated only if it vanishes. This is true, because Lyubeznik ([8], [9])
proved
$$\injdim (H^l_I(R))\leq \dim (\Supp_R (H^l_I(R)))$$ for any ideal $I$ and
any $l$ .Now if $0\neq H^l_I(R)$ was finitely generated, we would have
from [10], Theorem 18.9
$$\dim (R)=\depth (R)=\injdim (H^l_I(R))\leq \dim (\Supp_R (H^l_I(R)))\leq
\dim (R)$$ and consequently $\Supp_R (H^l_I(R))=\Spec (R)$ contradicting
$l\geq 1$.

In [3] Grothendieck conjectured that at least
$\Hom_R (R/I,H^l_I(R))$ is always finitely generated, but soon
Hartshorne was able to present the following counterexample to
Grothendieck's conjecture (see [4] for details and a proof): Let
$k$ be a field, $R=k[X,Y,Z,W]/(XY-ZW)=k[x,y,z,w]$, $I$ the ideal
$(x,z)\subseteq R$. Then $\Hom_R (R/I,H^2_I(R))$ is not finitely
generated.

However in Hartshorne's example the ring $R$ is not regular. Thus
the question arises whether Grothendieck's conjecture is true at
least in the regular case. In this context there is a theorem
([5], theorem 2.3(ii) and [8], corollary 3.5) stating that if $I$ is
an ideal of a regular ring $R$ which contains a field and $b$ is the
maximum of the heights of all primes minimal over $I$ then
for $l>b$, $\Hom_R (R/I,H^l_I(R))$ is finitely generated if and only if
$H^l_I(R)=0$.

Using this theorem one can give a counterexample to
Grothendieck's conjecture in the regular case, an idea which is
due to Hochster:

Let $k$ be a field of characteristic zero, $R=k[[X_1,\dots ,X_6]]$
a power series ring in six variables, $I_\Delta $ the ideal
generated by the $2\times 2$-minors of the matrix $\left( \matrix
{X_1 & X_2 & X_3 \cr X_4 & X_5 & X_6 \cr} \right) $. It can be
seen that $I_\Delta $ has pure height two and that $H^3_{I_\Delta
}(R)$ does not vanish. Now the above theorem implies
$\Hom_R (R/I,H^3_{I_\Delta }(R))$ is not finitely generated. But
theorem 7a) shows that at least the set of associated primes of
$\Hom_R (R/I,H^3_{I_\Delta }(R))$ (which is the same as
$\Ass_R (H^3_{I_\Delta }(R))$) is finite.

So one may wonder if any local cohomology module has only finitely
many associated primes. In [7] Huneke conjectured the following:
If $R$ is a local noetherian ring, then $\Ass_R (H^l_I(R))$ is
finite for any $I$ and any $l$. This paper deals with a weaker
version of Huneke's conjecture:
\bigskip
\mittel \noindent Conjecture ($*$): \normal
\medskip
\noindent If $R$ is a local Cohen-Macaulay-ring, then
$\Ass_R (H^l_I(R))$ is finite for any $I$ and any $l$.
\bigskip
\noindent Our main result is:
\bigskip
\mittel \noindent Theorem 6: \normal
\medskip
\noindent If $R$ is a local Cohen-Macaulay-ring, the following are
equivalent:

\noindent i) ($*$) is true for $R$.

\noindent ii) The following two conditions are fulfilled:

\noindent a) $\Ass_R (H^2_{(x,y)}(R))$ is finite for every $x,y\in R$.

\noindent b) $\Ass_R (H^3_{(x_1,x_2,y)}(R))$ is finite, whenever
$x_1,x_2\in R$ is a regular sequence and $y\in R$.
\bigskip
\noindent In Remark 2 it is shown that in the regular case
condition ii) a) is always satisfied. In fact at this point we will
not assume that $R$ is regular. We only need $R$ to be a so-called
almost factorial ring, which is weaker then being factorial.
\bigskip Besides this main result conjecture ($*$) is proved in
several special cases, for example in case $\dim (R)\leq3$ or furthermore
in case $\dim (R)\leq4$ provided $R$ is almost factorial.

Before going into the details, we remark that in the sequel we use
a certain (first-quadrant cohomological) spectral-sequence, the
socalled Groethendieck spectral-sequence for composed functors:

If $I$ and $J$ are ideals of a noetherian ring $R$, there is a
converging spectral-sequence
$$E^{p,q}_2=H^p_I(H^q_J(M))\Rightarrow H^{p+q}_{I+J}(M)$$ for
every $R$-module $M$: This is true because $\Gamma _J$ of an
injective module is injective again, where $\Gamma _J(M)$ is defined
as the submodule $\{ m\in M\vert J^n\cdot m=0\hbox{ for some }
n\} $ of $M$ (for details see [12], Theorem 5.8.3).
\parindent=0cm
\bigskip
We now start our examination of conjecture ($*$): At least for the spot
$l=\depth (I,R)$ there are only finitely many associated primes:
\bigskip
\mittel Theorem 1: \normal
\medskip
Let $(R,\goth m)$ be a noetherian local ring, $M$ a finitely
generated $R$-module and $I\subseteq R$ an ideal. Set
$t=\depth (I,M)$. Then $$\Ass_R (H^t_I(M))\subseteq
\Ass_R (\Ext^t_R (R/I,M))$$ and so $\Ass_R (H^t_I(M))$ is finite.
\medskip
Proof:

Choose $\goth p\in \Ass_R (H^t_I(M))$ arbritrarily. Because of
$H^t_{IR_\goth p}(M_\goth p)\neq 0$ we must have $t=\depth (IR_\goth
p,M_\goth p)$ and so we may assume $\goth p=\goth m$. Considering
the structure of $H^t_I(M)$ as a direct limit of certain
$\Ext $-modules we conclude $$\Hom_R (R/\goth m,\Ext^t_R (R/I^n,M))\neq
0$$ for some $n\in \Naturalsign $. Let $x_1,\dots ,x_t\in I$ be a
regular sequence. Using well-known formulas concerning $\Ext $ we
get $$\eqalign {0\neq \Hom_R (R/\goth m,
\Ext^t_R (R/I^n,M))&=\Hom_R (R/\goth m,\Hom_R (R/I^n,M/(x_1^n,\dots
,x_t^n)M))\cr &=\Hom_R (R/\goth m,\Hom_R (R/I,M/(x_1^n,\dots
,x_t^n)M))\cr &=\Hom_R (R/\goth m,\Ext^t_R (R/I,M))\ \ \ .\cr }$$
Now it follows that $\goth m\in \Ass_R (\Ext^t_R (R/I,M))$.\bigskip
A theorem established by M.P. Brodmann and A. Lashgari Faghani
([1], Proposition 2.1) states something more general: Let $R$
be a noetherian ring, $\goth a\subseteq R$ an ideal and $M$ a
finitely generated $R$-module. Furthermore, let $i\in \Naturalsign $
be given such that $H^j_{\goth a}(M)$ is finitely generated for
all $j<i$ and let $N\subseteq H^i_{\goth a}(M)$ be a finitely
generated submodule. Then, the set $\Ass_R (H^i_{\goth a}(M)/N)$ is finite.
\bigskip \mittel Lemma 1: \normal \medskip Let $R$ be a noetherian ring,
$M$ an $R$-module and $I,J$ ideals of $R$ with $\sqrt I\subseteq \sqrt
J$. Then $$H^l_I(M)=H^l_I(M/\Gamma _J(M))$$ for any $l\geq 1$.
\medskip
Proof:

Considering the long exact $\Gamma _I$-cohomology-sequence belonging to
$$0\longrightarrow \Gamma _J(M)\longrightarrow M\longrightarrow
M/\Gamma _J(M)\longrightarrow 0\ \ \ ,$$ we see it suffices to
show $H^l_I(\Gamma _J(M))=0$ for $l\geq 1$. Writing $M$ as the
union of its finitely generated submodules, we reduce to the case
$M$ itself is finitely generated, so that $\Gamma _J(M)$ is an
$R/J^n$-module ($n\gg 0$). Consequently $$H^l_I(\Gamma
_J(M))=H^l_{I(R/J^n)}(\Gamma _J(M))=H^l_{(0)}(\Gamma _J(M))=0\ \ \
.$$
\bigskip
Theorem 1 treated the case $l=\depth (I,R)$, and our next theorem
deals with the case $l=1$:
\bigskip \mittel
Theorem 2: \normal
\medskip
Let $R$ be a noetherian local ring, $I\subseteq R$ an ideal and
$M$ a finitely generated $R$-module. Then $\Ass_R (H^1_I(M))$ is
contained in $\Ass_R (\Ext^1_R (R/I,M/\Gamma _I(M)))$ and hence is
finite.
\medskip
Proof:

From Lemma 1 we get $$H^1_I(M)=H^1_I(M/\Gamma _I(M))$$ and $\Gamma
_I(M/\Gamma _I(M))=0$ implies $\depth (I,M/\Gamma _I(M))\geq 1$. So
theorem 2 becomes a corollary of theorem 1.
\bigskip
The next theorem shows that in studying conjecture ($*$), it
suffices to examine $H^j_I(R)$ when $\height (I)$ equals $j-1$
or $j$.
\bigskip
\mittel
Theorem 3:
\normal
\medskip
Let $(R,\goth m)$ be a local Cohen-Macaulay-ring, $I\subseteq R$
an ideal, $j>\height (I)$ and $H^j_I(R)\neq 0$. Then there exists an
ideal $\tilde I\supseteq I$ of height $j-1$ such that the natural
homomorphism $$H^j_{\tilde I}(R)\longrightarrow H^j_I(R)$$ becomes
an isomorphism.
\medskip
Proof:

We may assume $\height (I)<j-1$. Set $t=\height (I)$ and let
$x_1,\dots ,x_t\in I$ be a re\-gu\-lar sequence. We denote the
associated primes of $R/(x_1,\dots ,x_t)$ by $\goth p_1,\dots ,\goth
p_n$, enumerated in such a way that $$I\subseteq
\goth p_1\cap \dots \cap \goth p_r\ \ \ ,$$ $$I\nsubseteq \goth
p_{r+1},\dots ,\goth p_n\ \ \ .$$ We necessarily have $r<n$,
because otherwise $\sqrt I=\sqrt {(x_1,\dots ,x_t)}$ and
consequently $H^j_I(R)=0$, contrary to the assumptions. Using prime
avoidance we choose $$y\in (\goth p_{r+1}\cap \dots \cap \goth
p_n)\setminus (\goth p_1\cup \dots \cup \goth p_r)$$ and consider
the Mayer-Vietoris-sequence with respect to the ideals $(y),I$ and
the $R$-module $H^t_{(x_1,\dots ,x_t)}(R)=:M$: $$\eqalign
{&H^{j-t-1}_{I\cap (y)}(M)\longrightarrow
H^{j-t}_{(I,y)}(M)\longrightarrow
H^{j-t}_I(M)\oplus
H^{j-t}_{(y)}(M)\cr \longrightarrow
&H^{j-t}_{I\cap (y)}(M)\ \ \ .\cr }$$ In the sequel
we write $(\underline x)$ for the ideal $(x_1,\dots ,x_t)$ of $R$.
Because $j-t\geq 2$ and $I\cap (y)\subseteq \sqrt {(\underline
x)}$ it follows that $H^{j-t}_{(y)}=0$ and both the leftmost and
rightmost term in this sequence vanish; so the second arrow is an
isomorphism. Using the spectral-sequences for the composed
functors $\Gamma _{(I,y)}\circ \Gamma _{(\underline x)}$ and
$\Gamma _I\circ \Gamma _{(\underline x)}$ we conclude $$\eqalign
{H^j_{(I,y)}(R)&=H^{j-t}_{(I,y)}(M)\cr
&=H^{j-t}_I(M)\cr &=H^j_I(R)\ \ \ .\cr }$$
By construction $\height (I,y)=\height (I)+1$. Now the statement of
the theorem follows inductively.
\bigskip
The following corollary is the first step in a series of
reductions of conjecture ($*$): \bigskip \mittel Corollary 1:
\normal
\medskip
Let $(R,\goth m)$ be a local Cohen-Macaulay-ring and $j\in
\Naturalsign $. Then the following two statements are equivalent:

i) $\Ass_R (H^j_I(R))$ is finite for each ideal $I\subseteq R$.

ii) $\Ass_R (H^j_I(R))$ is finite for each ideal $I\subseteq R$
    satisfying $\height (I)=j-1$.\medskip
Proof:

Follows immediately from theorem 3.
\bigskip
Using the ideas of the proof of theorem 3 one can show that
$H^j_I(R)$ has only finitely many associated primes of
height $j$: \bigskip \mittel Corollary 2: \normal
\medskip
Let $(R,\goth m)$ be a local Cohen-Macaulay-ring, $I$ an ideal of
$R$ and $j\in \Naturalsign $. Then $$\Supp_R (H^j_I(R))\cap \{ \goth
p\in \Spec (R)\vert \height (\goth p)=j\} $$ is finite and therefore
$H^j_I(R)$ has only finitely many associated prime ideals of
height $j$.
\medskip
Proof:

We may assume $\height (I)\leq j-1$. Because of theorem 3 we may
even assume that the height of $I$ equals $j-1$. Let $x_1,\dots
,x_{j-1}\in I$ be a regular sequence and $\goth p_1,\dots ,\goth
p_n$ the associated primes of $R/(x_1,\dots ,x_{j-1})$, enumerated
in a way that we have $$I\subseteq \goth p_1\cap \dots \cap \goth
p_r\ \ \ ,$$ $$I\nsubseteq \goth p_{r+1},\dots ,\goth p_n\ \ \ .$$
We assume $r<n$ (if $r=n$ we have $\sqrt I=\sqrt {(\underline x)}$
and therefore $H^j_I(R)=0$). Set $J:=\goth p_{r+1}\cap \dots \cap
\goth p_n$ and consider the following part of a
Mayer-Vietoris-sequence: $$H^j_{I+J}(R)\longrightarrow
H^j_I(R)\oplus H^j_J(R)\longrightarrow H^j_{(x_1,\dots
,x_{j-1})}(R)=0\ \ \ .$$ It follows $\Supp_R (H^j_I(R))\subseteq
${\tensy \char "56}$(I+J)$. As $\height (I+J)\geq j$, there are only
finitely many primes of height $j$ in $\Supp_R (H^j_I(R))$.
\bigskip
The methods we have developed so far suffice to prove conjecture
($*$) in case $\dim (R)\leq 3$:
\bigskip
\mittel Corollary 3: \normal
\medskip Let $R$ be a local Cohen-Macaulay-ring of dimension at
most three, $I$ an ideal of $R$ and $j\in \Naturalsign $. Then
$H^j_I(R)$ has only finitely many associated primes.
\medskip
Proof:

Case $\dim (R)=2$: If $j=2$, the statement follows immediately from
theorems 1 and 3. The case $j=1$ is done by theorem 2.

Case $\dim (R)=3$: The case $j=3$ follows at once from theorems 1
and 3. $j=1$ is again done by theorem 2. If $j=2$, we assume
$\height (I)=1$ by theorem 2. Now the statement follows from
Corollary 2.
\bigskip
\mittel Lemma 2:
\medskip
\normal Let $I$ be an ideal of a noetherian ring $R$ and $M$ any
$R$-module. Then $\Ass_R (M/\Gamma _I(M))=\Ass_R (M)\cap
(\Spec (R)\setminus ${\tensy \char "56}$(I))$.
\medskip
Proof:

If $\goth p$ is associated to $M/\Gamma _I(M)$ we get from an
exact sequence $$0\longrightarrow R/\goth p\longrightarrow
M/\Gamma _I(M)$$ an exact sequence $$0\longrightarrow \Gamma
_I(R/\goth p)\longrightarrow \Gamma _I(M/\Gamma _I(M))=0$$ and
consequently $\goth p$ does not contain $I$. Choose $m\in M$
satisfying $\Gamma _I(M):m=\goth p$. Localizing we conclude
$$0:{m\over 1}=\Gamma _{IR_\goth p}(M_\goth p):{m\over 1}=\goth
pR_\goth p\ \ \ .$$ From our assumptions it follows that ${m\over
1}\neq 0$ , because otherwise there would exist $s\in R\setminus
\goth p$ with $sm=0$, contradicting $\Gamma _I(M):m=\goth p$.
Hence $\goth pR_\goth p\in \Ass_{R_\goth p} (M_\goth p)$,
equivalently $\goth p\in \Ass_R (M)$.

On the other hand, if we choose $\goth p\in \Ass_R (M)\cap
(\Spec (R)\setminus ${\tensy \char "56}$(I))$, $\goth p$ cannot be
associated to $\Gamma _I(M)$ and consequently must be associated
to $M/\Gamma _I(M)$ (consider $0\rightarrow \Gamma
_I(M)\rightarrow M\rightarrow M/\Gamma _I(M)\rightarrow 0$ exact).
\bigskip
\mittel Lemma 3: \normal
\medskip
Let $I$ be an ideal of a local Cohen-Macaulay-ring $R$ and set
$l=\height (I)+1$. Let $\goth p_1,\dots ,\goth p_n$ be the elements
of $\{ \goth p\in \Spec (R)\vert \goth p$ minimal over $I$ and
$\height (\goth p)=\height (I)\} $. Set $I^{pure}:=\goth p_1\cap \dots
\cap \goth p_n$. Then finiteness of $\Ass_R (H^l_{I^{pure}}(R))$
implies finiteness of $\Ass_R (H^l_I(R))$.
\medskip
Proof:

Let $\goth q_1,\dots ,\goth q_m$ be the elements of $\{ \goth p\in
\Spec (R)\vert \goth p$ minimal over $I$ and $\height (\goth
p)>\height (I)\} $ (without restriction assume $m\geq 1$) and set
$I^{\prime \prime }:=\goth q_1\cap \dots \cap \goth q_m$. Then
$\sqrt I=I^{pure} \cap I^{\prime \prime }$. Consider the
Mayer-Vietoris-sequence $$H^l_{I^{pure} +I^{\prime \prime
}}(R)\longrightarrow H^l_{I^{pure} }(R)\oplus H^l_{I^{\prime
\prime }}(R)\longrightarrow H^l_I(R)\longrightarrow
H^{l+1}_{I^{pure} +I^{\prime \prime }}(R)\ \ \ .$$ As by
construction $\height (I^{pure} +I^{\prime \prime })\geq
\height (I)+2=l+1$, the leftmost term in this sequence vanishes and
the rightmost term has only finitely many associated primes.
Furthermore $\height (I^{\prime \prime })\geq \height (I)+1=l$ and so
$H^l_{I^{\prime \prime }}(R)$ has only finitely many associated
prime ideals. Now the statement of the lemma is obvious.
\bigskip
Now we are in a position to give the next reduction of conjecture
($*$), which roughly spoken says one may restrict to the case $j=\mu
(I)$ when examining $\Ass_R (H^j_I(R))$: \bigskip \mittel Theorem 4:
\normal
\medskip Let $(R,\goth m)$ be a local Cohen-Macaulay-ring and
$t\in \Naturalsign $. Then the following two statements are
equivalent:

i) $H^{t+1}_I(R)$ has only finitely many associated prime ideals
for each ideal $I$ of $R$.

ii) Whenever $x_1,\dots ,x_t\in R$ is a regular sequence and $y\in
R$, the module $H^{t+1}_{(x_1,\dots ,x_t,y)}(R)$ has only finitely many
associated prime ideals.
\medskip
Proof:

Assume condition ii) is satisfied and let $I$ be an arbitrary
ideal of $R$. We have to show $\Ass_R (H^{t+1}_I(R))$ is finite.
Using Corollary 1 we may assume $\height (I)=t$. Using Lemma 3 we
can even assume that all primes minimal over $I$ have height $t$.

Let $x_1,\dots ,x_t\in I$ be a regular sequence and denote the
primes minimal over $I$ by $\goth p_1,\dots ,\goth p_n$. These are
also minimal over $(x_1,\dots ,x_t)$. Let $\goth q_1,\dots ,\goth
q_m$ be the other primes minimal over $(x_1,\dots ,x_t)$ (that is,
the ones not containing $I$). As all the ideals $\goth p_i$
and $\goth q_j$ have height $t$, we may choose a $$y^\prime \in
(\goth p_1\cap \dots \cap \goth p_n)\setminus (\goth q_1\cup \dots \cup
\goth q_m)\ \ \ .$$ Now a suitable power $y$ of $y^\prime $ will satisfy
$$y\in I\setminus (\goth q_1\cup \dots \cup \goth q_m)\ \ \ .$$ By
using Lemma 2 it follows that $y$ is not in any prime ideal
as\-socia\-ted to the $R$-mo\-du\-le $(R/(x_1^s,\dots
,x_t^s))/\Gamma _I(R/(x_1^s,\dots ,x_t^s))$ ($s\in \Naturalsign $
ar\-bi\-tra\-ry). Consequently $y$ operates injectively on
$(R/(x_1^s,\dots ,x_t^s))/\Gamma _I(R/(x_1^s,\dots ,x_t^s))$. From
the exactness of the direct limit-functor we conclude, that $y$
operates injectively on $$\eqalign {&\vtop{\baselineskip=1pt
\lineskiplimit=0pt \lineskip=1pt\hbox{$lim$}
\hbox{$\longrightarrow $} \hbox{$^{^{s\in \Naturalsign }}$}
}[(R/(x_1^s,\dots ,x_t^s))/\Gamma _I(R/(x_1^s,\dots ,x_t^s))]\cr
=&\vtop{\baselineskip=1pt \lineskiplimit=0pt
\lineskip=1pt\hbox{$lim$} \hbox{$\longrightarrow $}
\hbox{$^{^{s\in \Naturalsign }}$} }(R/(x_1^s,\dots ,x_t^s))/\Gamma
_I(\vtop{\baselineskip=1pt \lineskiplimit=0pt
\lineskip=1pt\hbox{$lim$} \hbox{$\longrightarrow $}
\hbox{$^{^{s\in \Naturalsign }}$} }(R/(x_1^s,\dots ,x_t^s)))\cr
=&H^t_{(x_1,\dots ,x_t)}(R)/\Gamma _I(H^t_{(x_1,\dots ,x_t)}(R))\
\ \ .\cr }$$ Call this property of $y$ ($**$). From well-known
spectral-sequence-arguments it follows $$\eqalign
{H^{t+1}_I(R)&=H^1_I(H^t_{(x_1,\dots ,x_t)}(R))\cr &\buildrel
(+)\over =H^1_I(H^t_{(x_1,\dots ,x_t)}(R)/\Gamma
_I(H^t_{(x_1,\dots ,x_t)}(R)))\cr &\buildrel (**)\over =\Gamma
_I(H^1_{(y)}(H^t_{(x_1,\dots ,x_t)}(R)/\Gamma _I(H^t_{(x_1,\dots
,x_t)}(R))))\cr &\subseteq H^1_{(y)}(H^t_{(x_1,\dots
,x_t)}(R)/\Gamma _I(H^t_{(x_1,\dots ,x_t)}(R)))\cr &\buildrel
(+)\over =H^1_{(y)}(H^t_{(x_1,\dots ,x_t)}(R))\cr
&=H^{t+1}_{(x_1,\dots ,x_t,y)}(R)\ \ \ .\cr }$$ The two equalities
(+) follow from Lemma 1. The above inclusion finishes our proof,
since we can conclude $$\vert \Ass_R (H^{t+1}_I(R))\vert \leq \vert
H^{t+1}_{(x_1,\dots ,x_t,y)}(R)\vert <\infty \ \ \ .$$
\bigskip
Using the various statements established so far, we can prove
conjecture ($*$) in the case $R$ is regular of dimension at most
four (cf. Theorem 5); in fact we do not actually need that $R$ is
regular. We will only use the fact that every height one prime
ideal is principal up to radical; this is true if $R$ is a Krull
domain whose divisor class group is torsion (cf. [2], Proposition
6.8). Krull domains whose divisor class group is torsion are usually
called almost factorial. In particular if $R$ is factorial, it is
almost factorial.
\bigskip
\mittel Theorem 5: \normal
\medskip
Let $R$ be a local almost factorial Cohen-Macaulay-ring of
dimension at most four, $I$ an ideal of $R$ and $j\in \Naturalsign
$. Then $H^j_I(R)$ has only finitely many associated primes, that
is, in these cases conjecture ($*$) is true.
\medskip
Proof:

We may restrict ourselves to the case $\dim (R)=4$. The case $j=0$
is trivial, $j=1$ follows from theorem 2, $j=3$ follows from our
corollaries 1 and 2 and $j=4$ from theorem 3. In the remaining
case $j=2$ we may assume $\height (I)=1$ (theorem 3). Using Lemma 3,
we may even assume that all primes minimal over $I$ have height
one. In our case this means that $I$ is principal up to radical
and so $H^2_I(R)=0$.
\bigskip
Theorem 6 is our final reduction of conjecture ($*$), allowing us to
restrict ourselves to the examination of "two" special cases (for
the regular case, see remark 2): \bigskip \mittel Theorem 6:
\normal
\medskip Let $R$ be a local Cohen-Macaulay-ring. Then the
following two statements are equivalent:

i) $H^j_I(R)$ has only finitely many associated prime ideals
for each ideal $I$ of $R$ and each $j\in \Naturalsign $.

ii) The following two conditions are satisfied:

a) $\Ass_R (H^2_{(x,y)}(R))$ is finite for every $x,y\in R$.

b) $\Ass_R (H^3_{(x_1,x_2,y)}(R))$ is finite whenever
$x_1,x_2\in R$ is a regular sequence and $y\in R$.
\medskip
Proof:

We only have to show ii) implies i): We do this by induction on
$j$:

$j=0$: Easy.

$j=1$: Theorem 2.

$j=2,3$: Theorem 4.

$j\geq 4$: Using theorem 4 we assume that $I=(x_1,\dots ,x_j)$
(for some $x_1,\dots ,x_j\in R$). Here $[\ ]$ means Gaussian
brackets, that is $[q]:=max\{ i\in {\bf Z} \vert i\leq q\} $ for rational
$q$. Set $I^\prime :=(x_1,\dots,x_{[j/2]}),
I^{\prime \prime }:=(x_{[j/2]+1},\dots ,x_j)\subseteq
R$ ideals and consider the following Mayer-Vietoris-sequence:
$$H^{j-1}_{I^\prime }(R)\oplus H^{j-1}_{I^{\prime \prime
}}(R)\longrightarrow H^{j-1}_{I^\prime \cap I^{\prime \prime
}}(R)\longrightarrow H^j_I(R)\longrightarrow H^j_{I^\prime }\oplus
H^j_{I^{\prime \prime }}(R)\ \ \ .$$ Combined with our induction
hypothesis (using $j-1\geq j-([j/2]+1)+1$) we get from this an
isomorphism $$H^{j-1}_{I^\prime \cap I^{\prime \prime
}}(R)\longrightarrow H^j_I(R)\ \ \ .$$ Another application of our
induction hypothesis finishes the proof of the theorem.
\bigskip
\mittel Remark 1: \normal
\medskip
i) Let $R$ be a local Cohen-Macaulay-ring, $n\in \{ 2,3\} $ and
$x_1,\dots ,x_n\in R$. Now from $\vert \Ass_R (H^n_{(x_1,\dots
,x_n)}(R))\vert <\infty $ conjecture ($*$) would follow. We can
write the module $H^n_{(x_1,\dots ,x_n)}(R)$ in another way. First
we have $$H^n_{(x_1,\dots
,x_n)}(R)=H^1_{(x_1)}(H^{n-1}_{(x_2,\dots ,x_n)}(R))$$ and from
the right-exactness of $H^1_{(x_1)}$ we may conclude
$$H^1_{(x_1)}(H^{n-1}_{(x_2,\dots ,x_n)}(R))=H^1_{(x_1)}(R)\otimes
_RH^{n-1}_{(x_2,\dots ,x_n)}(R)\ \ \ .$$ An easy induction proof
gives us $$H^n_{(x_1,\dots ,x_n)}(R)=H^1_{(x_1)}(R)\otimes _R\dots
\otimes _RH^1_{(x_n)}(R)=(R_{x_1}/R)\otimes _R\dots \otimes
_R(R_{x_n}/R)\ \ \ .$$ So for conjecture ($*$) it is sufficient to
prove $$\vert \Ass_R ((R_{x_1}/R)\otimes _R\dots \otimes
_R(R_{x_n}/R))\vert <\infty $$ for $n\in \{ 2,3\} $.

ii) Consider the complete case, that is, $R$ is a local complete
Cohen-Macaulay-ring. Similar to theorem 6, condition ii) assume
$t\in \{ 1,2\} $, $x_1,\dots ,x_t\in R$ a regular sequence and $y\in R$.
Consider $R$ as an $R[[T]]$-module via the $R$-algebra homomorphism
$R[[T]]\longrightarrow R$ sending $T$ to $y$. We then calculate
$$\eqalign {H^{t+1}_{(x_1,\dots ,x_t,y)}(R)&=H^{t+1}_{(x_1,\dots
,x_t,T)}(R)\cr &=H^{t+1}_{(x_1,\dots ,x_t,T)}(R[[T]]/(T-y))\cr
&=H^{t+1}_{(x_1,\dots ,x_t,T)}(R[[T]])/(T-y)H^{t+1}_{(x_1,\dots
,x_t,T)}(R[[T]])\ \ \ .\cr }$$ Since $x_1,\dots ,x_t,T\in R[[T]]$
is a regular sequence, it is in the complete case sufficient (for
conjecture ($*$)) to show that whenever $t\in \{ 2,3\} $, $x_1,\dots
,x_t\in R$ is a regular sequence and $y\in R$ we have $$\vert
\Ass_R (H^t_{(x_1,\dots ,x_t)}(R)/yH^t_{(x_1,\dots ,x_t)}(R))\vert
<\infty \ \ \ .$$
\bigskip
\mittel Remark 2: \normal
\medskip
If $R$ is an almost factorial local ring, condition a) from
theorem 6 ii) is automatically fulfilled. To show this we may, with
respect to theorem 3, assume $\height (x,y)=1$. Using Lemma 3 we may
even assume that all primes minimal over $(x,y)$ have height one.
As $R$ is almost factorial, it follows that $(x,y)$ is principal up to
radical and so $H^2_{(x,y)}(R)=0$.
\bigskip
The remaining theorems 7 and 8 prove conjecture ($*$) in certain
generic cases (where $R/I$ is Cohen-Macaulay); theorem 7 treats
the equicharacteristic case and theorem 8 deals with mixed
characteristics.
\bigskip \mittel Theorem 7: \normal
\medskip
a) let $k$ be a field, $R=k[[X_1,\dots ,X_6]]$ a power series ring
in six indeterminates, $\Delta _1:=X_2X_6-X_3X_5,\Delta
_2:=X_1X_6-X_3X_4,\Delta _3:=X_1X_5-X_2X_4$ (these are the
$2\times 2$-minors of the matrix $\left( \matrix {X_1 & X_2 & X_3
\cr X_4 & X_5 & X_6 \cr} \right) $), $I$ the ideal $(\Delta
_1,\Delta _2,\Delta _3)\subseteq R$. Then
$\Supp_R (H^3_I(R))\subseteq \{ (X_1,\dots ,X_6)\} $ and
consequently $\Ass_R (H^3_I(R))$ is finite.

b) Let $R$ be a local equicharacteristic Cohen-Macaulay-ring and
$x_1,\dots ,x_6\in R$ be a regular sequence. Let $\delta
_1:=x_2x_6-x_3x_5,\delta _2:=x_1x_6-x_3x_4,\delta
_3:=x_1x_5-x_2x_4$ and $I$ be the ideal $(\delta _1,\delta _2,\delta _3)
\subseteq R$. Then $\Ass_R (H^3_I(R))$ is finite.
\medskip
Proof:

a) It is well-known that $R/I$ is a Cohen-Macaulay domain of
dimension 4. Consequently $I$ is a prime ideal of height two. From
[10], Theorem 30.4.(ii) it follows that
$$Sing(R/(\Delta _1))\subseteq \{ \goth p\in \Spec (R/(\Delta
_1))\vert \goth p\supseteq (X_2,X_6,X_3,X_5)\} \ \ \ .$$ Here
$Sing(R/(\Delta _1))$ means the set of all primes $\goth p$
satisfying $(R/(\Delta _1))_\goth p$ is not regular. Furthermore
we have $$Sing(R/(\Delta _2))\subseteq \{ \goth p\in
\Spec (R/(\Delta _1))\vert \goth p\supseteq (X_1,X_6,X_3,X_4)\} $$
and $$Sing(R/(\Delta _3))\subseteq \{ \goth p\in \Spec (R/(\Delta
_1))\vert \goth p\supseteq (X_1,X_5,X_2,X_4)\} \ \ \ .$$ Choose
$\goth p\in \Spec (R/I)\setminus \{ (X_1,\dots ,X_6)\} $
arbitrarily. We have to show $H^3_{IR_\goth p}(R_\goth p)=0$. From
our above calculations we know there is an $i\in \{ 1,2,3\} $ with
$\goth p\notin Sing(R/(\Delta _i))$. Thus $(R/(\Delta _i))_\goth
p$ is factorial. Combining this with the fact that $I/(\Delta _i)$
is a prime ideal of height one, we conclude the
ideal $IR_\goth p/(\Delta _i)R_\goth p\subseteq R_\goth p/(\Delta
_i)R_\goth p$ is principal. This finally shows
$$0=H^2_{IR_\goth p/(\Delta _i)R_\goth p}(H^1_{(\Delta _i)R_\goth
p}(R_\goth p))=H^3_{IR_\goth p}(R_\goth p)\ \ \ .$$ b) We may
assume that $R$ is complete, because if the statement is proved in
the complete case, then the formula $$\Ass_R (H^3_I(R))=\bigcup
_{\goth p\in \Ass_R (H^3_IR(R))}\Ass_{\hat R} (\hat R/\goth p\hat R)$$
(cf. [10], Theorem 23.2.(ii)) implies finiteness of
$\Ass_R (H^3_I(R))$ (each $\Ass_{\hat R} (\hat R/\goth p\hat R)$
contains a $\goth q$ with $\goth q\cap R=\goth p$).

Let $k\subseteq R$ be a field, $k[[X_1,\dots ,X_6]]$ be a power
series ring in six variables and $\Delta _1,\Delta _2,\Delta _3\in
k[[X_1,\dots ,X_6]]$ (like in a)) the $2\times 2$-minors of $
\left( \matrix {X_1 & X_2 & X_3 \cr X_4 & X_5 & X_6 \cr} \right)$.
The flat $k$-algebrahomomorphism $$k[[X_1,\dots
,X_6]]\longrightarrow R$$ with $X_i\mapsto x_i$ ($i=1,\dots ,6$)
sends $\Delta _j$ to $\delta _j$ ($j=1,2,3$). This implies
$$H^3_I(R)=H^3_{(\Delta _1,\Delta _2,\Delta _3)}(R)=H^3_{(\Delta
_1,\Delta _2,\Delta _3)}(k[[X_1,\dots ,X_6]])\otimes
_{k[[X_1,\dots ,X_6]]}R$$ and we conclude
$$\Ass_R (H^3_I(R))\subseteq \Ass_R (R/(X_1,\dots ,X_6)R)\ \ \ ,$$
from [10], Theorem 23.2.(ii), which finally proves b).
\bigskip
\mittel Theorem 8: \normal
\medskip
a) Let $p$ be a prime number, $C$ a complete $p$-ring,
$R=C[[X_1,\dots ,X_6]]$ a power series ring in six variables and
set $\Delta _1:=X_2X_6-X_3X_5,\Delta _2:=X_1X_6-X_3X_4,\Delta
_3:=X_1X_5-X_2X_4$ (these are the $2\times 2$-minors of the matrix
$\left( \matrix {X_1 & X_2 & X_3 \cr X_4 & X_5 & X_6 \cr} \right)
$), $I$ the ideal $(\Delta _1,\Delta _2,\Delta _3)\subseteq R$.
Then $\Supp_R (H^3_I(R))\subseteq ${\tensy \char "56}$((X_1,\dots
,X_6))$ and consequently $\Ass_R (H^3_I(R))$ is finite.

b) Let $p$ be a prime number, $(R,\goth m)$ be a local
Cohen-Macaulay-ring satisfying $char(R)=0$, $char(R/\goth m)=p$
and $x_1,\dots ,x_6\in R$ with the property that $p,x_1,\dots
,x_6\in R$ is a regular sequence. Set $\delta
_1:=x_2x_6-x_3x_5,\delta _2:=x_1x_6-x_3x_4,\delta
_3:=x_1x_5-x_2x_4$ and let $I$ be the ideal $(\delta _1,\delta _2,
\delta _3) \subseteq R$. Then $\Ass_R (H^3_I(R))$ is finite.
\medskip
Proof:

a) The proof is practically the same as the proof of theorem 7 a).

b) Like in the proof of theorem 7 b), we may assume that $R$ is
complete. According to [10], theorem 29.3 $R$ has a coefficient
ring $C\subseteq R$. Let $C[[X_1,\dots ,X_6]]$ be a power series
ring in six variables and $\Delta _1,\Delta _2,\Delta _3\in
C[[X_1,\dots ,X_6]]$ (like in a)) the $2\times 2$-minors of
$\left( \matrix {X_1 & X_2 & X_3 \cr X_4 & X_5 & X_6 \cr} \right)
$. The rest of the proof may be copied from the proof of theorem 7
b) until one finally gets $$\Ass_R (H^3_I(R))\subseteq
\Ass_R (R/(X_1,\dots ,X_6)R)\cup \Ass_R (R/(p,X_1,\dots ,X_6)R)\ \ \
,$$ which proves b).

\def\litem{\par\noindent \hangindent=\parindent\ltextindent}
\def\ltextindent#1{\hbox to \hangindent{#1\hss}\ignorespaces}
\vfill \eject \pageno=15 \gross References\normal
\bigskip \bigskip
\parindent=0.8cm
\litem{1.} Brodmann, M.P. and Lashgari Faghani, A. A finiteness result
for associated primes
of local cohomology modules, {\it preprint}, (1998).
\litem{2.} Fossum, R. M. The Divisor Class Group of a Krull
Domain, {\it Springer-Verlag}, (1973).
\medskip
\litem{3.} Grothendieck, A. Cohomologie locale des faiscaux
coh\'{e}rents et th\'{e}or\`{e}mes de Lefschetz locaux et globaux,
{\it S.G.A.} {\bf II}, (1968).
\medskip
\litem{4.} Hartshorne, R. Affine duality and cofiniteness, {\it
Inventiones Mathematicae} {\bf 9}, (1970), 145-164.
\medskip
\litem{5.} Huneke, C. and Koh, J. Cofiniteness and vanishing of
local cohomology modules, {\it Math. Proc. Camb. Phil. Soc.} {\bf
110}, (1991), 421-429.
\medskip
\litem{6.} Huneke, C. and Sharp, R. Bass Numbers of Local
Cohomology Modules, {\it Transactions of the American Mathematical
Society} {\bf 339}, (1993).
\medskip
\litem{7.} Huneke, C. Problems on Local Cohomology, {\it Res.
Notes Math. } {\bf 2}, (1992).
\medskip
\litem{8.} Lyubeznik, G. Finiteness properties of local
cohomology modules (an application of $D$-modules to Commutative
Algebra), {it Inventiones Mathematicae} {113}, (1993).
\medskip
\litem{9.} Lyubeznik, G. $F$-Modules: Applications to Local
Cohomology and $D$-modules in Characteristic $p>0$, preprint.
\medskip
\litem{10.} Matsumura, H. Commutative ring theory, {\it Cambridge
University Press}, (1986).
\medskip
\litem{11.} Singh, A. $p$-torsion elements in local cohomology
modules, preprint, (1999).
\medskip
\litem{12.} Weibel, C. A. An introduction to homological
algebra, {\it Cambridge University Press}, (1994).

\end